\documentclass[12pt]{article}
\usepackage[english]{babel}
\usepackage{amsmath,cite,amsthm}
\usepackage{amssymb}
\usepackage{latexsym}
\newtheoremstyle{theorem}%name
  {10pt}          % space above
  {10pt}  % space below
  {\sl}  % bofy font
  {\parindent}     % ident - empty=no indent,  \parindent= paragraph indent
  {\bf}  % thm head font
  {. }    % punctuation after thm head
  { }    % space after thm head: `` ``=normal \newline=linebreak
  {}     % thm head specification
\theoremstyle{theorem}
\newtheorem{theorem}{Theorem}

\newtheoremstyle{defi}%name
  {10pt}          % space above
  {10pt}  % space below
  {\rm}  % bofy font
  {\parindent}     % ident - empty=no indent,  \parindent= paragraph indent
  {\bf}  % thm head font
  {. }    % punctuation after thm head
  { }    % space after thm head: `` ``=normal \newline=linebreak
  {}     % thm head specification
\theoremstyle{defi}

\begin{document}

\title{$\mathbb{B}_{0}$-VALUED MONOGENIC FUNCTIONS\\ %AND\\ THEIR APPLICATIONS 
TO THE THEORY OF\\ PLANE ANISOTROPY} %\\ ANISOTROPIC PLANE MEDIA}

\author{S.\,V.~Gryshchuk\\
 Institute of Mathematics, \\ National Academy of Sciences
of Ukraine,\\ Tereshchenkivska Str. 3, 01004, Kyiv, Ukraine\\
gryshchuk@imath.kiev.ua, serhii.gryshchuk@gmail.com}

\maketitle

\begin{abstract}
A solution of the elliptic type PDE of the 4th order, being  a reduction of the Eqs. of stress function
corresponding to any case of plane anisotropy which is not equal to isotropy (proved by S.\,G.~Mikhlin),  
is described in terms of hypercomplex ``analytic'' functions with values in two-dimensional semisimple algebra over the field of complex numbers in case when a domain under consideration is bounded and simply-connected. A boundary value problem 
 on finding a function  which satisfies this PDE %the stress Eqs.
  in the considered domain (bounded and simply-connected) and permits continuations (to the boundary) of operators $\frac{\partial}{\partial x}$ and $\frac{\partial}{\partial y}$ acted to it,
 is reduced to certain BVP for these hypercomplex functions.
 % problem (have a great importance in the 2nd main BVP of Anisotropy by S.\,G.~Mikhlin)

 \vspace*{6mm}
 
 %Primary 30G35, 31A30; Secondary
 %74B05 %2010

 % 74B05 Classical linear elasticity,  31A30 Biharmonic, polyharmonic functions and equations, Poisson's equation
%30G35 Functions of hypercomplex variables and generalized variables

\end{abstract}

\newpage
%%%%%%%%%%%%%%%%%%%%%%%%%SECTION-1
\section{Introduction}

The algebraic-analytic approaches to the investigation of elastic media in terms of ``analytic'' functions satisfying
a system of partial differential equations (a generalization of the ``Cauchy--Riemann conditions'') with
values in finite-dimensional algebras were developed in \cite{Sodbero,Edenhofer,Gilb-Wendl-Proc-75,Kov-Mel,Kov,GrPl_umz-09,Gr-Pl_Dop2009,Cont-13,Gr-Amade-15,
 mon-f-bih-BVP-MMAS,zb17Gr,openMath17G,Gr-dop-15,Bon} (commutative algebras and isotropic plane media),
\cite{p-bikh-al,umzhOrth18-1,umzhOrth18-2} (commutative algebras and orthotropic two-dimensional media), \cite{gurl-ed-el+qu-poth-14,Gurl-Kol-Mushf-MMA15,Grigorev-2017,Tsalik86,Tsalik95,Gurl-Spr-bo-Appl-hol-16} (algebra of quaternions and
isotropic three-dimensional media),
\cite{Gurl-Spr-bo-Appl-hol-16,Liu-Cliff-3D-Elast} (approaches of using different kinds of monogenic functions with values in the Clifford algebras for solving the equilibrium system of the
isotropic three-dimensional media),
 \cite{Sold-hyperan-SMiPr-04,Sold-SMFN-dr_nar-16,Abop-Sold-Lame_orthopr-SMiPr-04} (algebras of complex (2$\times$2) matrices and anisotropic plane media),
and \cite{Mitin-98} (algebras of complex (3$\times$3) matrices and anisotropic plane media).

The present paper is devoted to the construction of classes of ``analytic'' functions $\Phi$ with values in two-dimensional commutative algebras over the field of complex numbers containing bases $(e_{1}, e_{2})$ with some algebraic properties
(in what follows, we construct all mentioned bases and the corresponding algebra in the explicit form) sufficient for the real components of these functions to satisfy the following equations for fixed $p > 0$,
$p \ne 1$:
%\[ %
\begin{equation} \label{genBihEqMikhl}
\widetilde{l}_{p} u(x,y):= \left(\frac{\partial^4 }{\partial y^4}+
A_{p} \frac{\partial^4
}{\partial x^2\partial y^2}+
B_{p}\frac{\partial^4}{\partial x^4}
 \right) u(x,y)=0, %\qquad e_1^2+e_2^2\ne 0,
 %\eqno (1) \]%
\end{equation}
where $A_{p}:= p^2 + 1$, $B_{p} = p^2$, $u$ is a real-valued solution of
\eqref{genBihEqMikhl}, an argument $(x,y) \in D$, while the latter is  belonging to the Cartesian plane $xOy$.

The operator $\widetilde{l}_{p}$  can be factorizated in the form:
\begin{equation}\label{factOPorthMikhl}
  \widetilde{l}_{p} = \widetilde{l}_{1,p} \circ \Delta_2 = \Delta_2 \circ \widetilde{l}_{1,p},\,\,
  \widetilde{l}_{1,p}:= \frac{\partial^2}{\partial y^2} + p^2 \frac{\partial^2}{\partial x^2},
\end{equation}
where $\widetilde{l}_{1,p} \circ \Delta_2$ is a symbol of composition of operators $\widetilde{l}_{1,p}$ and $\Delta_2$,
 $\Delta_2:= \frac{\partial^2}{\partial x^2}+
\frac{\partial^2}{\partial y^2}$ is the 2-D Laplasian.

Equation \eqref{genBihEqMikhl} is a special case of the generalized biharmonic equation (this term was used, e.g., in \cite{Mikhlin} and \cite[p.~603]{Mush_upr}),
which is extremely important in the anisotropic two-dimensional theory of elasticity (see \cite{Mikhlin-anisotr,Mush_upr,Fridman-50,Lekh-TU-an-tel,Sherman-ani-Se-38,
 Bog-izv-05,Parton-meth-math-th-upr,kupradze-potentsial}) and determines (in absence of body forces)
the equation for finding the stress function $u(x, y)$~(in the isotropic case a similar function is often called the Airy function and Eqs.~\eqref{genBihEqMikhl} turns into the biharmonic equation).  %with p = 1

S.\,G.~Mikhlin proved in \cite{Mikhlin-anisotr} that an equation of finding the stress function in any cases of plane anisotropy (except of isotropic case)
can be reduced to  Eqs.~\eqref{genBihEqMikhl}. From the other side, formally Eqs. \eqref{genBihEqMikhl}
 corresponds to orthotropic material being an equation of finding the stress function in
a special case of plane anisotropy~--- orthotropy (cf., e.g., \cite[pp.~33,34]{Lekh-TU-an-tel}, \cite{Gr-Lodz-18}).

Note that the Eqs.~\eqref{genBihEqMikhl} is considered  with $A_{p}:=2p$, $B_{p}:=1$, $p>1$ in \cite{umzhOrth18-1,umzhOrth18-2},
or $-1<p<1$ in \cite{ProcIPMM18}.
%%%%%%%%%%%%%%%%%%%%%%%%%%%%%%%%%%%%%%%%%%%%%%%%%%%%%%%%%%%%%%%%%%%%%%%%%%%%%%%%%SECTION-2
%%%%%%%%%%%%%%%%%%%%%%%%%%%%%%%%%%%%%%%%%%%%%%%%%%%%%%%%%%%%%%%%%%%%%%%%%%%%%%%%SECTION-1
\section{Two-dimensional algebras over the field of complex numbers and their bases
associated with Eqs.~(1.1)}
It is known (see \cite{Study}) that there exist (to within isomorphism) two associative algebras of the second rank
with identity e commuting over the field of complex numbers $\mathbb{C}$ . These algebras are generated by the bases $(e,\rho)$ and
 $(e, \omega)$, respectively:
\begin{equation}\label{alg-B}
\mathbb{B}:=\{c_1 e+c_2\rho: c_k\in\mathbb{C},k=1,2\},\,\rho^2=0,
\end{equation}
\begin{equation}\label{alg-Bo}
\mathbb{B}_{0}:= \{c_1 e+c_2 \omega: c_k\in\mathbb{C},k=1,2\},\,
\omega^2=e.
\end{equation}
It is clear that the algebra $\mathbb{B}_{0}$  is semisimple (for the definition, see, e.g., \cite[p.~33]{Cheb}) and contains a basis with
orthogonal idempotents $(\mathcal{I}_1, \mathcal{I}_2)$, where
\begin{equation}\label{bas-idempot}
\mathcal{I}_1=\frac{1}{2}\left(e+\omega\right),\,
\mathcal{I}_2=\frac{1}{2}\left(e-\omega\right),\,
\mathcal{I}_1\mathcal{I}_2=0.
\end{equation}
It is clear that
\begin{equation}\label{idemp-e-omeg}
\mathcal{I}_1+\mathcal{I}_2=e,\,
\mathcal{I}_1-\mathcal{I}_2=\omega.
\end{equation}
In the works of different %foreign
researchers, several names are used for the algebra $\mathbb{B}_{0}$.
 Thus, in \cite{Sobch-unipodal-alg}, it is called
{\it unipodal}. Moreover, it determines the simplest case of complex Clifford algebra (cf., e.g.,  \cite{Sobch-unipodal-alg,bayl-cliff-alg}).
Algebra (3) is a complexification of the algebra of hyperbolic or double numbers $P$ over the field of real
numbers $\mathbb{R}$:
\[\mathbb{B}_{0}=
\mathbb{P}\oplus i \mathbb{P},\, \mathbb{P}:=\{x e+h y:
x,y\in\mathbb{R}\}, h:=\omega,\]
where $i$ is the imaginary unit.

The element $w = c_{1}\mathcal{I}_{1} + c_{2}\mathcal{I}_{2}$ from $\mathbb{B}_{0}$
 is invertible if and only if $c_{k} \ne  0, k = 1, 2$. In this case, the inverse
element is given by the equality (see \cite[p.~38]{bayl-cliff-alg})
\begin{equation}\label{B0-inver}
w^{-1}=\frac{1}{c_1}\,\mathcal{I}_1+
\frac{1}{c_2}\,\mathcal{I}_2.
\end{equation}
Since the algebra $\mathbb{B}$ contains a nonzero radical $\{c\rho: c \in \mathbb{C}\}$ (see \cite{Kov-Mel}), the algebra $\{c\rho: c \in \mathbb{C}\}$ is not semisimple.
An element $a = c_{1}e + c_{2} \rho$ from $\mathbb{B}$ is invertible if and only if $c_1\ne 0$.
 In this case, the equality
$a^{-1}=\frac{1}{c_1}\,e-\frac{c_2}{(c_1)^2}\,\rho$
is true (see \cite{Gr-Pl-Zb-10}).

%Let p > 1 be an arbitrary fixed real number, p 2 R.
For any complex number $s$ we introduce the
notation
\begin{equation}\label{lPtilde}
  \widehat{l}_{p}(s):= s^4 + \left(p^2 + 1\right) s^2 + p^2.
\end{equation}

The equation $ \widehat{l}_{p}(s)=0$ is the characteristic equation of the
\eqref{genBihEqMikhl}, its set of roots is
\begin{equation}\label{rootsChMikh}
  \{s_1, s_2, \overline{s_1}, \overline{s_2}\} =:  \ker \widehat{l}_{p}, \, s_1 = i, s_2 = ip,
\end{equation}
where $\overline{x+iy}:= x-iy$, $x, y \in \mathbb{R}$.

Now we are looking for an associative, commutative algebra of the second rank with unity $e$  over the field of complex numbers $\mathbb{C}$ and containing a basis $(e_{1}, e_{2})$ that satisfies the condition
\begin{equation} \label{tab_umn_ba-p-orthMikhl}
\mathcal{L}_{p}(e_1,e_2): = e_{2}^4+ G_{p}\,
e_{1}^2 e_{2}^{2} +  H_{p}\, e_{1}^4 = 0.
\end{equation}
where
$G_{p}:=\left(p^2 + 1\right)$, $H_{p}:=p^2$.

Note that a similar problem had been considered  for the
equation of the type like \eqref{tab_umn_ba-p-orthMikhl} with $G_{p}:= 2p$, $H_{p}:=1$,
 $p>1$ in \cite{umzhOrth18-1} and $-1<p<1$ in \cite{ProcIPMM18}.

Doing in analogous way as in the proof of similar Theorem in  \cite{umzhOrth18-1} one can obtain the following theorem.

%\vspace{3mm} \noindent \textbf{Theorem 1.} \textit{
\begin{theorem}\label{baseOrtMikh}
The algebra $\mathbb{B}$ does not contain any basis $(e_{1}, e_{2})$ satisfying condition \eqref{tab_umn_ba-p-orthMikhl}.
There exists  a set of cardinality continuum of
 bases $(e_{1}, e_{2})$ in $\mathbb{B}_0$ satisfying condition \eqref{tab_umn_ba-p-orthMikhl}\emph{:}
 \begin{equation}\label{counInfSetOrtMikh}
   e_1=\alpha \mathcal{I}_{1}+ \beta \mathcal{I}_{2},\,
e_{2}=\alpha \widetilde{s}_1\,\mathcal{I}_{1}+ \beta \widetilde{s}_2\,\mathcal{I}_{2} \,\,
\forall \alpha, \beta \in \mathbb{C} \setminus\{0\},
 \end{equation}
where  $\widetilde{s}_{k} \in \ker \widetilde{l}_{p}$, $k=1,2$,
$\widetilde{s}_{1} \ne \widetilde{s}_{2}$ %}
\end{theorem}

Let us restrict our attention on the case
 $\alpha = \beta \equiv 1$, $\widetilde{s}_{k}={s}_{k}$, $k=1,2$, in \eqref{counInfSetOrtMikh}. Therefore,
 we have
 \begin{equation}\label{Case-base-orthMihl}
   e_1 = \mathcal{I}_{1} +\mathcal{I}_{2}  \equiv e,\,
   e_2 = i\left( \mathcal{I}_{1} + p \mathcal{I}_{2}
   \right)  \equiv  i\left(e_1 + (p-1) \mathcal{I}_{2}
   \right).
 \end{equation}

 %It follows from \eqref{Case-base-orthMihl}
 Since expressions of idempotents $\mathcal{I}_{k}$, $k=1, 2$, via elements
 of bases \eqref{Case-base-orthMihl} are
 \begin{equation}\label{idempOrtMikh}
 \mathcal{I}_{1} = - \frac{1}{1-p}\left(p e_1 + i e_2 \right),\,
\mathcal{I}_{2} = \frac{1}{1-p}\left(e_1 + i e_2 \right),
 \end{equation}
we obtain the multiplication table for the bases
\eqref{Case-base-orthMihl}:
\begin{equation}\label{tablMult-basOrtMikh}
  e_1 e_2 = e_2, \, %\, e_1 =e, \,
   e_2^2 = p\, e_1 + i(p+1)\, e_{2}.
\end{equation}
%%%%%%%%%%%%%%%%%%%%%%%%%%%%%%%%%%%%%%%%%%%%%%%%%%%%%%%%%%%%%%%%%%%%%%%%%%%%%%%%%%%%%%%%%SECTION-3
%%%%%%%%%%%%%%%%%%%%%%%%%%%%%%%%%%%%%%%%%%%%%%%%%%%%%%%%%%%%%%%%%%%%%%%%%%%%%%%%%%%%%%%%%

\section{$\mathbb{B}_{0}$-valued monogenic functions and Eqs.~{\bf (1)} }%\eqref{genBihEqMikhl}
Consider  $\mu_{e_1,e_2}=\{\zeta=x\,e_1+y\,e_2 :
x,y\in\mathbb R\}$ which is a linear span of the elements $e_1,e_2$
of the basis (\ref{Case-base-orthMihl}) over the field of real
numbers $\mathbb R$. With a domain $D$ of the Cartesian plane $xOy$
we associate the congruent domain $D_{\zeta}:= \{\zeta=xe_1+ye_2 :
(x,y)\in D\} \subset \mu_{e_1,e_2}$, and corresponding
domains in the complex plane $\mathbb{C}$: $D_{z}:=\{z=x+iy:(x,y)\in
D\}$, $D_{z_{p}}:=\{z_{p}=x+ipy:(x,y)\in
D\}$. Let $D_{\ast}$ be a domain in $xOy$ or in $\mu_{e_1,e_2}$.
 Denote by $\partial D_{\ast}$
 a boundary of a
domain $D_{\ast}$, $\mathrm{cl}D_{\ast}$ means a closure of a domain
$D_{\ast}$.

%Denote by $\partial D$ ($\partial D_{\zeta}$) a boundary of a
%domain $D$ ($D_{\zeta}$), $\mathrm{cl}D$ ($\mathrm{cl}D_{\zeta}$)
%means a closure of a domain $D$ ($D_{\zeta}$).

In what follows, $(x,y)\in D$, $\zeta=x\,e_1+y\,e_2\in
D_{\zeta}$, $z=x+iy\in D_{z}$, $z_{p}=x+ipy \in D_{z_{p}}$.

Inasmuch as divisors of zero don't belong to $\mu_{e_1,e_2}$,
one can define the derivative $\Phi'(\zeta)$ of function $\Phi
\colon D_{\zeta}\longrightarrow \mathbb{B}_{0}$ in the same way as in
the complex plane:
\[\Phi'(\zeta):=\lim\limits_{h\to 0,\, h\in\mu_{e_1,e_2}}
\bigl(\Phi(\zeta+h)-\Phi(\zeta)\bigr)\,h^{-1}\,.\]
 We say that a function $\Phi \colon D_{\zeta}\longrightarrow \mathbb{B}$ is
\textit{monogenic} in a domain $D_{\zeta}$ and, denote as
$\Phi\in\mathcal{M}_{\mathbb{B}_{0}}(D_{\zeta})$, if the derivative
$\Phi'(\zeta)$ exists in every point $\zeta\in D_{\zeta}$.

Every function  $\Phi \colon D_{\zeta}\longrightarrow \mathbb{B}_{0}$
has a form
\begin{equation}\label{mon-funk} \Phi(\zeta)=U_{1}(x,y)\,e_1+U_{2}(x,y)\,ie_1+
U_{3}(x,y)\,e_2+U_{4}(x,y)\,ie_2,
\end{equation}
where $\zeta=xe_1+ye_2$, $U_{k}\colon D\longrightarrow \mathbb{R}$,
$k=\overline{1,4}$.

Let denote every real component  $U_{k}$, $k=\overline{1,4}$, in expansion
 (\ref{mon-funk})  by
$\mathrm{U}_{k}\left[\Phi\right]$, i.~e., for arbitrary fixed
$k\in\{1,\dots,4\}$:
$$\mathrm{U}_{k}\left[\Phi(\zeta)\right]:=U_{k}(x,y)\,\, \forall \zeta= xe_1 + y e_2 \in D_{\zeta}.$$

 We establish the following theorem similar to analogous theorems in \cite{GrPl_umz-09,umzhOrth18-1}.

%\vspace{3mm} \noindent \textbf{Theorem~2.}
 %\textit{
 \begin{theorem}\label{CR-OrtMikh}
  A function $\Phi \colon D_{\zeta} \longrightarrow \mathbb{B}_{0}$ is monogenic in the $D_{\zeta}$ if and only if its components
$U_{k} \colon D \longrightarrow \mathbb{R}, k = \overline{1,4}$, in decomposition \eqref{mon-funk} are differentiable in the domain D and the following analog
of the Cauchy--Riemann conditions is true\emph{:}
\begin{equation}\label{usl_K_R}
\frac{\partial \Phi(\zeta)}{\partial y}\,e_1=\frac{\partial
\Phi(\zeta)}{\partial x}\,e_2.
\end{equation}
\end{theorem}
%}

In an extended form the condition (\ref{usl_K_R}) for the monogenic
function (\ref{mon-funk}) is equivalent to the system of four
equations (cf., e.g., \cite{Kov-Mel,GrPl_umz-09}) with respect to
components $U_{k}=\mathrm{U}_{k}\left[\Phi\right]$,
$k=\overline{1,4}$, in (\ref{mon-funk}):
\begin{eqnarray}
\frac{\partial U_{1}(x,y)}{\partial y}&=& p \, \frac{\partial
U_{3}(x,y)}{\partial x},\label{kr1orMikh}\\ %\nonumber
\frac{\partial U_2(x,y)}{\partial y}&=& p \, \frac{\partial U_{4}(x,y)}{\partial x},\label{kr2orMikh}\\
%{}\label{K-R_Meln}\\
 \frac{\partial U_{3}(x,y)}{\partial y}&=&\frac{\partial U_{1}(x,y)}{\partial x}-\left(p+1\right)\frac{\partial U_{4}(x,y)}{\partial x},\label{kr3orMikh}\\
\frac{\partial U_{4}(x,y)}{\partial y}&=&\frac{\partial
U_{2}(x,y)}{\partial x}+\left(p+1\right)\frac{\partial U_{3}(x,y)}{\partial
x}\label{kr4orMikh}.
\end{eqnarray}

Using \eqref{Case-base-orthMihl}, an element $\zeta = x e_1 + y e_2 \in \mu_{e_1,e_2}$ turns of by the formula
\begin{equation}\label{zetaI-12}
  \zeta = z \, \mathcal{I}_{1} + z_{p}\, \mathcal{I}_{2} \, \forall \zeta \in \mu_{e_1,e_2}.
\end{equation}

A function $\Phi\in\mathcal{M}_{\mathbb{B}_{0}}(D_{\zeta})$ can be expressed in terms of two holomorphic functions of the complex variable $z$ and $z_{p}$, respectively.
The following theorem obtained with use of \eqref{zetaI-12}  similar to analogous theorem in \cite{umzhOrth18-1}.

%\vspace{3mm} \noindent \textbf{Theorem~3.}
 %\textit{
 \begin{theorem}\label{predMonOrtMikh}%{
 The function $\Phi \colon D_{\zeta} \longrightarrow \mathbb{B}_{0}$ is monogenic in the domain $D_{\zeta}$ if and only if the following
equality is true:
  \begin{equation}\label{rep-mon-f-B-0-orMikh}
\Phi(\zeta)=
F_{1}\left(z\right)\mathcal{I}_{1}+
F_{2}\left(z_{p}\right)\mathcal{I}_{2}\,\, \forall \zeta
\in D_{\zeta},
\end{equation}
where $F_{1}(z)$, $F_{2}(z_p)$  is a holomorphic function of its complex variable $z\in D_{z}$, $z_{p}\in D_{z_{p}}$, respectively.%}
\end{theorem}

It follows from \eqref{rep-mon-f-B-0-orMikh} and Theorem~2  %\eqref{usl_K_R}
that every derivative $\Phi^{(n)}$, $n=1,2,\dots$, is monogenic.
Thus, we have $\widetilde{l}_{p}\Phi(\zeta) =\mathcal{L}_{p}(e_1,e_2)
 \Phi^{(4)}(\zeta) \equiv 0 $.
Therefore, we deduce that every component $U_{k}=\mathrm{U}_{k}\left[\Phi\right]$, $k=\overline{1,4}$, satisfies Eqs.~\eqref{genBihEqMikhl}.

 With use of  \eqref{idempOrtMikh} we rewrite  \eqref{rep-mon-f-B-0-orMikh} to the form
 \[ \Phi(\zeta)=\frac{1}{1-p}\left(\left(-p \, F_{1}(z) + F_{2}(z_{p})\right) e_1 +
 \left(F_{2}(z_{p})-\frac{1}{p} \, F_{2}(z)\right)\right) \, \forall \zeta \in D_{\zeta}. \]
After that, substituting with loss of generality $c_{k}F_{k}$ to $F_{k}$, $k=1,2$, where
$c_{1}:= -\frac{p}{1-p}$, $c_{2}:= \frac{1}{1-p}$, we obtain the following representation of the monogenic function  $\Phi$ in the bases \eqref{Case-base-orthMihl} for every $\zeta \in D_{\zeta}$:
\begin{equation}\label{repMonE1E2OrtMikh}
  \Phi(\zeta) = \left(F_{1}(z) + F_{2}\left(z_{p}\right)\right)e_{1}+
  i \left(F_{2}\left(z_{p}\right)  -  \frac{1}{p}\, F_{1}\left(z\right)\right)e_{2}.
\end{equation}

Since now we assume that $D$ is a bounded and simply-connected domain.

Than  by solving a system \eqref{kr1orMikh}~--~\eqref{kr4orMikh} with $U_{1}\equiv 0$ it is easy to deliver that
a function $\Phi_{1,0}\in\mathcal{M}_{\mathbb{B}_{0}}(D_{\zeta})$, such
that $\mathrm{U}_{1}\left[\Phi_{1,0}\right] \equiv 0$, has a form
\begin{equation}\label{phi1-0OrtMikh}
\Phi_{1,0}(\zeta) = a i\left(y \, e_1 + \left(\frac{z_{p}}{p}+ \frac{i}{p}\, y \right)\, e_2 \right) +
b i e_1 + c e_2 + d i e_2 \, \, \forall \zeta \in D_{\zeta},
\end{equation}
 where $a, b, c, d$ are arbitrary real numbers.

 We shall prove that for every fixed solution $u$ of the equation \eqref{genBihEqMikhl} in a bounded simply
connected domain $D\subset xOy$   exists a function $\Phi_{u} \in\mathcal{M}_{\mathbb{B}_{0}}(D_{\zeta})$  such that $\mathrm{U}_{1}\left[\Phi_{u}\right] \equiv u$.

There is a well-known fact  (cf., e.g., \cite[\S 20, p.~136]{Lekh-TU-an-tel} or \cite{Sherman-ani-Se-38}),
  that every solution $u$ %the general solution
  of the equation \eqref{genBihEqMikhl}
is expressed in the form:
\begin{equation}\label{pr-bikhP-2AnalitFun-sOrtMikh}
  u(x,y) = \mathrm{Re\,}\left(F_{1}\left(z\right)+F_{2}\left(z_{p}\right)\right)\,\,
  \forall (x,y) \in D,
\end{equation}
where
$F_{1} \colon D_{z} \longrightarrow \mathbb{C}$
and $F_{2} \colon D_{z_{p}} \longrightarrow \mathbb{C}$
are  analytic functions of their variables. %any

By use of
\eqref{repMonE1E2OrtMikh} with $\Phi_{u}:=\Phi$ and $F_{k}$ the same as in \eqref{pr-bikhP-2AnalitFun-sOrtMikh}, $k=1,2$, we rewrite the equality
 \eqref{pr-bikhP-2AnalitFun-sOrtMikh} in the form
 \begin{equation}\label{reprGenSolVia1ComponOrMikh}
   u(x,y)= \mathrm{U}_{1}\left[\Phi_{u}(\zeta) \right] \, \forall \zeta \in D_{\zeta},
 \end{equation}
 where  $\Phi_{u} \in \mathcal{M}_{\mathbb{B}_{0}}(D_{\zeta})$.

It follows now from
\eqref{phi1-0OrtMikh} and
\eqref{pr-bikhP-2AnalitFun-sOrtMikh} the following theorem
being an analog of the classical fact that any  harmonic function (in the bounded simply-connected domain of the real plane) is a real part of an
analytic function of the complex variable, moreover, this representation is unique up to the imaginary constant as an
 addend.

%\vspace{3mm} \noindent \textbf{Theorem~4.}
 %\textit{
 \begin{theorem}\label{1CompOrtMikh}
Let $u$ be a solution of the equation~\eqref{genBihEqMikhl}.
Then all $\Phi \in \mathcal{M}_{\mathbb{B}_{0}}(D_{\zeta})$,
such that
\begin{equation}\label{U1PhiuOrMikh}
  u(x,y)= \mathrm{U}_{1}\left[\Phi(\zeta) \right] \, \forall \zeta \in D_{\zeta},
\end{equation}
are expressed by the formula
\begin{equation}\label{sumOrMikh}
  \Phi(\zeta) = \Phi_{u}(\zeta) + \Phi_{1,0}(\zeta) \, \forall \zeta \in D_{\zeta}.
\end{equation}%}
\end{theorem}

Note, that  similar Theorem is proved in
 with deal of the equation of the type like~\eqref{genBihEqMikhl}
with $A_{p}:=2p$, $B_{p}:=1$ and $p>1$ or $-1<p<1$ (being  Eqs.~of the stress function to a certain class of orthotropic
plane deformations) in \cite{umzhOrth18-2}($p>1)$ or \cite{ProcIPMM18}($-1<p<1$).

%%%%%%%%%%%%%%%%%%%%%%%%%%%%%%%%%%%%%%%%%%%%%%%%%%%%%%%%%%%%%%%%%%%%%%%%%%%%%%%%SECTION-3
\section{A BVP of plane anisotropy and corresponding  BVP for $\mathbb{B}_{0}$-valued monogenic functions}

Consider a boundary value problem on finding a function $u \colon D \longrightarrow \mathbb{R}$
satisfying conditions
\begin{equation}
    \label{genBikhBVP}
\left\{
\begin{array}{ll}
 \widetilde{l}_{p} u(x,y)= 0 \, \, \forall (x,y) \in D, \vspace*{2mm} % & \forall z \in \partial D,
\\
 \lim\limits_{(x,y) \to (x_{\hspace*{-0.4mm}\circ},y_{\hspace*{-0.4mm}\circ})\in \partial D, (x,y) \in D}\frac{\partial u(x,y)}{\partial x} = u_{1}(x_{\hspace*{-0.4mm}\circ},y_{\hspace*{-0.4mm}\circ})\, \, \forall (x_{\hspace*{-0.4mm}\circ},y_{\hspace*{-0.4mm}\circ}) \in \partial D,
 \\
 \lim\limits_{(x,y) \to (x_{\hspace*{-0.4mm}\circ},y_{\hspace*{-0.4mm}\circ})\in\, \partial D, (x,y) \in D}\frac{\partial u(x,y)}{\partial y} = u_{3}(x_{\hspace*{-0.4mm}\circ},y_{\hspace*{-0.4mm}\circ})\, \, \forall (x_{\hspace*{-0.4mm}\circ},y_{\hspace*{-0.4mm}\circ}) \in \partial D,
\end{array}
\right.
\end{equation}
where $u_{k} \colon \partial D \longrightarrow \mathbb{B}_{0}$, $k\in\{1,3\}$,
are given continuous  functions.

The problem \eqref{genBikhBVP} has a great importance in the anisotropy theory (cf., e.g., \cite{Mikh_int_eq_Th-El}), when $\widetilde{l}_{p}$ is a biharmonic operator we arrive at the isotropic case and this changed boundary value problem \eqref{genBikhBVP} named as the {\it biharmonic problem} (cf., e.g., \cite{Mikhlin})].  There are different approaches to its solving (see \cite{Mikhlin-anisotr,Sherman-ani-Se-38,bogan-11,Mikh_int_eq_Th-El}).

Our aim is to find a new method of solving BVP~\eqref{genBikhBVP} by use of  monogenic functions.
Let $\Phi_1 \in \mathcal{M}_{\mathbb{B}_{0}}(D_{\zeta})$ such that
\begin{equation}\label{1Comp}
 \mathrm{U}_{1}[\Phi_{1}(\zeta)] = u(x,y) \, \, \forall (x,y) \in D,
\end{equation}
where $u$ is a south-fought solution of Problem~\eqref{genBikhBVP}.  Using the
equality  \eqref{kr1orMikh} and the condition \eqref{usl_K_R}, we have
\[\Phi_{1}'(\zeta) \equiv \frac{\partial \Phi_1(\zeta)}{\partial x} =
\frac{\partial u(x,y)}{\partial x}\, e_1 + \frac{1}{p}\frac{\partial u(x,y)}{\partial x}\, e_2 +\frac{\partial \mathrm{U}_{2}\left[\Phi_{1}(\zeta)\right]}{\partial x} ie_1 +
 \frac{\partial \mathrm{U}_{4}\left[\Phi_{1}(\zeta)\right]}{\partial x} ie_2 \]
for all $\zeta \in D_{\zeta}$, therefore, for all $(x,y)\in D$.
Thus, BVP~\eqref{genBikhBVP} is reduced to BVP on finding a monogenic function $\Phi:=\Phi'$
satisfying boundary conditions
\[\lim\limits_{\zeta \to \zeta_{\circ}= x_{\hspace*{-0.25mm}\circ}\, e_{1} + y_{\hspace*{-0.25mm}\circ}\,e_{2} \in \partial D_{\zeta}, \zeta \in D_{\zeta}} \mathrm{U}_{k}[\Phi(\zeta)] = \lambda_{k} u_{k}(x_{\hspace*{-0.4mm}\circ},y_{\hspace*{-0.4mm}\circ}), k\in\{1,3\},\, \, \forall (x_{\hspace*{-0.4mm}\circ},y_{\hspace*{-0.4mm}\circ}) \in \partial D, \]
where $\lambda_{k}:= 1$ if $k=1$, $\lambda_{k}:= \frac{1}{\,p\,}$ if $k=3$.
Solving the latter BVP for monogenic finctions, we obtain a solution of BVP~\eqref{genBikhBVP} in the form
\[
u(x,y)=\int_{(x_{\hspace*{-0.4mm}\circ},y_{\hspace*{-0.4mm}\circ})}^{(x,y)}\left(
\mathrm{U}_{1}[\Phi(xe_1+ye_2)]\,dx\right. + \left. p \,\mathrm{U}_{3}[\Phi(xe_1+ye_2)]\,dy\right) + \mathrm{const}_{\mathbb{R}} \,\, \forall (x,y) \in D,
\]
where
$\mathrm{const}_{\mathbb{R}}$ is an arbitrary real number,
$(x_{\hspace*{-0.4mm}\circ},y_{\hspace*{-0.4mm}\circ})$
is a fixed point in $D$,
integration means along any piecewise
smooth curve jointing %, which joints
this point with a point with variable coordinates
$(x, y)$.

\vspace{3mm}

\section*{ACKNOWLEDGEMENTS}
This research is partially supported by the State Program of Ukraine (Project\break
No.~0117U004077).
%%%%%%%%%%%%

%\newpage

\end{document}